# Smoothness and Smooth Extensions (I): Generalization of MWK Functions and Gradually Varied Functions

Li Chen
Department of Computer Science and Information Technology
University of the District of Columbia
lchen@udc.edu
May 19, 2010

**Abstract**

A mathematical smooth function means that the function has continuous derivatives to a certain degree $C^{(k)}$. We call it a $k$-smooth function or a smooth function if $k$ can grow infinitively. Based on quantum physics, there is no such smooth surface in the real world on a very small scale. However, we do have a concept of smooth surfaces in practice since we always compare whether one surface is smoother than another one. This paper deals with the possible definitions of "natural" smoothness and their relationship to the original mathematical definition of smooth functions. The motivation of giving the definition of a smooth function is to study smooth extensions for practical applications.

We observe this problem from two directions: From discrete to continuous, we suggest considering both micro smooth, the refinement of a smoothed function, and macro smooth, the best approximation using existing discrete space. (For two-dimensional or higher dimensional cases, we can use Hessian matrices.) From continuous to discrete, we suggest a new definition of natural smooth, it uses a scan from down scaling to up scaling to obtain the a ratio for sign changes by ignoring zero to represent the smoothness. For differentiable functions, mathematical smoothness does not mean a "good looking" smooth for a sampled set in discrete space. Finally, we discuss the Lipschitz continuity for defining the smoothness, which will be called discrete smoothness. This paper gives philosophical consideration of smoothness for practical problems, rather than a mathematical deduction or reduction, even though our inferences are based on solid mathematics.

## 1. Introduction

The world can be described by continuity and discontinuity. Most parts are continuous. Discontinuity provides the boundary of continuous components, and continuity gives the content of discontinuous. An image can be separated as several continuous components. The edge of a component indicates the discontinuity to two parties. However, the boundary itself could be continuous. Continuity sometimes means controllable and discontinuity means uncontrollable such as something suddenly happened. For instance, a sunny day suddenly had a rain.

Continuous sometime is not good enough to describe a phenomenal. Therefore, in continuous, we need to define smoothness. As well as in discontinuous, we want to define chaos. This paper only focuses on the concept of smooth and its extensions.

Mathematicians always try to make things simpler. What is the smoothness? The mathematical smooth function means that a function has continuous derivatives to a certain degree $C^{(k)}$. We call it as a $k$-smooth function or the smooth function if $k$ could be infinitively growth. Based on quantum physics, there is no such smooth surface in real world in quantum-scale. So the "natural" smoothness is relative not absolute.

On the other hand, a mathematical smooth function does not always mean the continuity in real world. For instance, when a sunny day starts to rain. The rain usually does not come by a centimeter in a minute. The rain usually had few drops then become bigger. So in a microscope, the weather changes continuously or smoothly per say. But the truth is that the sunny day had a rain! That is discontinued.

It comes a problem, what shall we define the smoothness that can be used to real world problems and it does not course a contradiction to our beautiful mathematics exist. In fact, mathematician has already observed the problem. They have used Lipschitz conditions [15] and rectifiability to restrict a function to have a better continuous and smooth [11].

It was perfect amazing that three papers published in the same year dealing with a same or very similar problem. McShane, Whitney, and Kirszbraun all did independently in 1934. They obtained an important theorem for the Lipschitz function extension. What they proved is: For any subset $E$ of a metric space $S$, if there is function $f$ on $E$ satisfies the Lipschitz condition, then it can be extended to $S$ preserving the same Lipschitz constant.

So our problem is solved! Actually it is not. McShane-Whitney-Kirszbraun Theorem only provides a theoretical result. This is because that the Lipschitz constant could be very big. No one who is involved in real problems wants to admit such a function is a "continuous" function.

In 1983, Krantz published a paper that specifically studied smoothness of functions in Lipschitz space for technical applications [14]. Unfortunately, such an important movement did not get enough attentions in mathematics.

Along with the fast development of computer science, scientists started to study "continuous functions" in a digital format for image processing. Rosenfeld pioneered this direction. In1982, he used a fuzzy logic concept to study the relationship among pixels. He defined the concept of fuzzy connected components. In 1985, Chen proposed another type of connectedness between two pixels called lambda-connectedness. This concept was found later having closer connections to continuous functions. In 1986, Rosenfeld proposed discretely continuous functions [16]. In 1989, Chen simplified lambda-connectedness to gradual variation and proved the necessary and sufficient condition of the extension [6]. In 1989, Steele used the Lipschitz condition to check the image "smoothness."

Chen's constructive method suggests an algorithm that is totally different from McShane's construction that builds a function from a maximum or minimum possibilities. The algorithm designed by Chen uses the local and middle construction that is more suitable to the practice. In addition, Chen used a sequence of real or rational numbers instead of integers, treating the problem more generally. That overcomes the limitation of Lipschitz space by enlarging to local Lipschitz [1-8].

What is our point? Even though we have a continuous or mathematical smooth function to describe the addressed question, the problem may only be able to be presented in a finite form for large amount of applications using finite samples in today's technology.

Based on McShane-Whitney-Kirszbraun theorem, a continuous surface always exists. According to the classic the Weierstrass Approximation Theorem: If $f$ is a continuous real-valued function on *[a,b]*, a closed and bounded interval, $f$ can be uniformly approximated on that interval by polynomials to any degree of accuracy. In addition to the results obtained by Fefferman [12], differentiable extensions of the Whitney's problem are exists. However, the solutions are based on infinitively refinement of the domain space. In practice, we have the limited availability for the space we are dealing with. Recall the example of the sunny day had a rain, even though we could insert the time to micro second, we still have a discontinuous event that the sunny day usually does not have a rain.

The mathematical definition of smooth functions does not provide a precise description of such phenomenal of real problems. The existence of a smooth function for "sudden rain" does not good enough to the problem of our sampling scale by minutes, or second. Smoothness is depending on sampling in real world or nature! For many cases absolute smoothness based on the mathematical definition may not be working right.

This paper deals with the possible definitions of natural smoothness and their relationship with the original mathematical definition of smooth functions. The motivation of giving the definitions of smooth function is to study smooth extension for practical applications.

We observe this problem in two directions: From discrete to continuous, we suggest to consider both of micro smooth meaning the refinement to a smoothed function, and macro smooth meaning the best approximation using existing discrete space. For two-dimensional or high dimensional cases, we can use the Hessian matrix to determine the extreme points and to represent the case.

From continuous to discrete, we suggest a new definition of natural smooth, it uses a scan from down scaling to up scaling to get the a ratio for sign changes with ignoring zero to represent the smoothness. For differentiable functions, mathematical smoothness does not mean a "good looking" smooth for a sampled set in discrete space. Finally, we discuss the Lipschitz continuity for defining the smoothness. This paper is to give philosophical consideration of smoothness for practical problem, rather a mathematical deduction or reduction even though our inferences are based on solid mathematics.

Since 2005, Fefferman started to publish papers for Whitney's smooth extensions and suggest the technology to solve the real problems of finite sample points. This inspired our recent research. Our resent research results suggest that we need to go one more step to clarify what real smooth function is in real world.

## 2. Review of Concepts

We have reviewed many existing methods for continuous functions' extension [2]. Despite many important practical methods, we first introduce here the key theorem by McShane, Whitney, and Kirszbraun. It gives a construction to a continuous extension to a set including a finite set. McShane gave a constructive proof for the existence of the extension in [15]. He constructed a minimal extension (INF) that is Lipschitz.

**Theorem 2.1** (The McShane-Whitney-Kirszbraun Theorem) Let J be a closed subset of D, and f(J) be a function. If $|f(x) - f(y)| \leq Lip|x - y|$ where $Lip$ is a constant for all x, y in J, then

$$F(x) = \max_{a \in J} \{f(a) - Lip|a - x|\} \tag{2.1}$$

is a continuous extension such that $|F(x) - F(y)| \leq Lip|x - y|$ for all x, y in D.

We can also construct a maximum extension by

$$F(x) = \min_{a \in J}\{f(a) + Lip|a-x|\} \tag{2.2}$$

We can call it a maximum extension (SUP). It is obvious that neither INF nor SUP can be directly used in data reconstruction. However, *F=(INF+SUP)/2* can be a reasonable function and *F* is a Lipschitz extension. We will call *F=(INF+SUP)/2* the McShane-Whitney-Kirszbraun (MWK) *mid function* [4]. All these functions used the similar method to construct is called MWK functions.

In 1989, Chen obtained a more general result for the extension on a finite set [6]. In fact, in 1986, Rosenfeld had realized the need of such an extension for image processing [16]. Chen defined the concept of gradual variation is a discrete method that can be built on any graph. The gradually varied surface is a special discrete surface. We now introduce this concept.

Let function *f: D→{$A_1, A_2,...,A_n$}*, if *a* and *b* are adjacent in *D* implies *f(a)=f(b)*, or *f(b) =$A_{(i-1)}$* or *$A_{(i+1)}$* when *f(a)=$A_i$* , point *(a,f(a))* and *(b,f(b))* are said to be gradually varied. A 2D function (surface) is said to be gradually varied if every adjacent pair are gradually varied. We can see that the gradually varied function is more general than a Lipschitz function.

Discrete Surface Fitting: Given *J⊆D*, and *f: J→{$A_1, A_2,...,A_n$}* decide if there is a *F: D→{$A_1, A_2,...,A_n$}* such that *F* is gradually varied where *f(x)=F(x)*, *x* in *J*.

**Theorem 2.2** (Chen, 1989[6]) The necessary and sufficient conditions for the existence of a gradually varied extension *F* is: for all *x,y* in *J*, *d(x,y)≥ |i-j|*, *f(x)=$A_i$* and *f(y)=$A_j$*, where d is the distance between *x* and *y* in *D*.

The construction method selected by Chen was different from McShane's paper[15] or McShane-Whitney-Kirszbraun mid function. Chen's method was to expand points by making no change or making the minimum change the base value. It maintains good properties naturally. Examples are shown below [2]. The MWK mid function may different from gradually varied functions, see examples in [2].

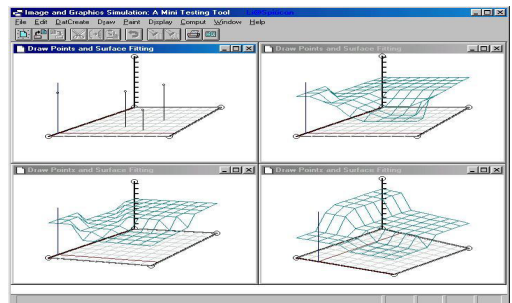

Fig 1. Gradually varied interpolations: guiding points and the gradually varied fitting.

These two theorems presented above guarantee the existence of the continuous extension. But to get the smooth extension is a problem. Based on the discussion in Section 1, the mathematical smooth extension for a continuous function exists in terms of approximation. However, such an extension does not provide a looking smooth function since human's eye does not have such a good micro smooth extension. Our eye can only sample certain resolution (digitization) in a scale.

# 3. Micro- and Macro- Smoothness: From Discrete functions to Continuous Functions

Based on Whitney's original work in 1930s and recent Fefferman's important accomplishments [12], there always be a smooth function that is the extension for giving finite sample points. The question is that, in applications, we could not use refinement infinitively. If we do not want to change the resolution/scale or to change it very limited to the original resolution, how do we get a smooth (looking) function? We must to re-define some definitions including smoothness that relates to error estimation for uniformed approximation.

In industry, there already some considerations exist which are called micro smoothness and macro smoothness. They are used for measuring the similar things in real world problems. In paper and fabric related industry, Dimmick described that *macro-smoothness* is based on base paper deformation while *micro-smoothness* is based on coating [10]. Micro is a Local smoothness, and macro is a global smoothness. Micro is coated smooth while macro is material and intrinsic smooth. Micro is external smooth and macro is internal smooth [10][19].

We observe this problem in two directions: From discrete to continuous, we suggest to consider both of micro smooth meaning the refinement to a smoothed function, and macro smooth meaning the best approximation using existing discrete space.

Let us assume a discrete function $f: D \to R$. How do we make $F: D' \to R$ that is smoother than $f$ if $D$ is a subset of $D'$ and $|D'| < c\,|D|$? Where $c$ is a constant.

## 3.1 Micro Smoothness

If we try to get a micro smooth function, we could use Fefferman's method [12]. The method Fefferman adopted is to use Whitney's Jets at each point, then use a binary tree structure to combine them. Each step of the union operation is to generate a system of linear inequalities based on the Lipschitz conditions for every degree of derivatives. A Whitney's Jet is a local expansion of Taylor series at a refined grid or meshes.

## 3.2 Macro Smoothness

The macro smooth function is to keep (or make a little change) the original grid or meshes. If we assume the original function is Lipschitz with Lipschitz constant *Lip*, we have to make changes to get a smooth function in terms of using finite difference instead of derivatives. Use approximation (limited errors) instead of equal left and right derivatives.

Here is something to consider: when we make a macro smooth function, there is no need to change values for every point; we can select some points that will not be changed (similar to guiding points in application). In such a case, it is equivalent to the case of extension of finite sampling cases. We have used this philosophy in our algorithm design in [1][2][3].

The algorithm to complete this type of task may need iterations. The finite difference error is the key to consider. We will discuss its qualitative aspect later.

## 3.3 Relationship with existing concepts

In discussion of natural smoothness of continuous functions, the micro smoothness is the "smoothness" for high resolution (dense sampling) and the macro smoothness is the smoothness for low resolution. One could say these are the high frequency and the low frequency components. What is the difference?

We are not very interested in every single frequency like Fourier and wavelet analysis. We are more interested in the extreme cases such as maximum and minimum (possible "frequencies"). That are micro and macro "smoothness." But where to find them is a key. We could use persistent analysis to determine them if necessary. Thus,

$$Function = micro\ component + macro\ component \qquad (3.1)$$

This is an approximation to the real data.

## 4. Natural Smoothness of Continuous Functions --- From continuous functions to discrete functions

Since a continuous function has a uniform approximation from sequence of polynomial functions. The mathematical smoothness for a continuous function does not really mean much. From continuous to discrete, we suggest a new definition of natural smoothness, it uses a scan from large scaling to small scaling to get the ratio for sign changes of derivatives on the curve with ignoring zero to represent the smoothness. For differentiable functions, mathematical smoothness does not mean a "good looking" smooth for a sampled set in discrete space. If the function is a rectifiable curve, such a sampling always exists [11].

A function that is not continuous could be a Lipschitz function that could have a small Lipschitz constant. It may have a good natural smoothness. In other words, after sampling, the macro smoothness is good. The micro smoothness may not be good. It is just like that material smoothness (such as for paper and cloth) is related to deformable properties.

Another good example is that a mountain looks smooth on a skyline, when we get close to the mountain we see a lot of trees that is not smooth at all. So the natural smoothness is related to resolutions and scales. For one-dimensional function, assume that:

*Nsamples* denotes the number of samples selected in the function.

*DerivativeSignChanges* denotes the total sign changes of the derivative function of the reconstructed function based on the samples by the finite difference, the MWK mid function, or gradually varied functions.

**Definition 4.1:** Natural smoothness of a continuous *1D* function is a (stable) ratio:

$$R = (Nsamples - DerivativeSignChanges) / Nsamples \qquad (4.1)$$

The worst case is the function with the pattern at sampling points $\{-1,+1,-1,+1,\ldots\}$. So the smoothness is almost 0. Limits of this formula or the probability of the estimation would be the fact. The statistical method such as student test may be needed to find the confidence region. When a small set is sampled, the micro-changes will not be affected the value of the natural smoothness. This will overcome the inconsistence of classical mathematical definition of smooth functions. Natural smoothness is a concept with relativity.

Let's define the natural smoothness for 2D and high dimensional cases. We still use sampling and reconstruction and then to finding smoothness. There are two steps to get the smoothness: (1) Get sample points, *SN*, and their values. (2) Reconstruct the function using MWK mid function or the gradually varied function (The algorithm used in [1-3]). (3) Use Hessian matrices to find extreme points of the reconstructed functions. Count the number *EN*.

The sampling is necessary since we want to eliminate the high frequency part based on the average sampling scale. Another method of using traditional wavelets or Fourier transforms to get the decomposition of frequencies, then reconstruction based on different frequency scales can be made. We still can calculate the numbers of extreme points using Hessian matrices.

**Definition 4.2:** Natural smoothness of a continuous k-*D* function is a (stable) ratio:

$$R=(SN-EN)/EN \qquad (4.1)$$

## 5. Natural Smooth Functions for Discrete Sets and Gradually Varied Smooth Reconstruction

If we just have a finite numbers of samples and no other information, we may only be able to fit a continuous function using gradually varied fitting. If we want it to be a good smooth looking, we could polish the function. Based on the definition of natural smoothness, the linear interpolation (if one dimensional), or gradually varied fitting and polishing will not change the natural smoothness of the fitted function. They are essentially the similar functions.

Using gradually varied functions for smooth reconstruction can be either micro- or macro-smooth reconstruction. Since for most of applications, we like to keep the original grid, and most of existing reconstruction used polynomial method that will yield a mathematically smooth function, so we mainly use gradually varied functions for macro smooth functions.

We have discussed natural smoothness in Section 4. The natural smoothness would change when the grid scale changes if using gradually varied functions. What we want to obtain is that the fitted function has the best uniform approximation to the guiding points.

This method we have used recently [1][2][3]. The easiest way is to use the finite difference method to calculate the left and right derivatives for higher order of derivatives. This finite difference method will not change the extreme points. So it shall and will keep the natural smoothness. Natural smooth is a type of macro in down scaling, and is Micro in up (dense) scaling.

Since the case is in the discrete set, it is hard to get exact same value of the left derivative and right derivative. We only want to limit the difference between them for a given error epsilon. Related ideas are presented in [14][12]. We want,

$$| (f(x_{i+1})-f(x_i))/(x_{i+1}-x_i) - (f(x_{i+2})-f(x_{i+1}))/(x_{i+2}-x_{i+1})| \leq epsilon \qquad (5.1)$$

If f is Lipschitz, and $x_i$ are equally arranged, we have

$$| 2 f(x_{i+1}) - f(x_i) - f(x_{i+2}) | \leq epsilon$$

Krantz had estimation for an example in [14].

In summary, without additional information such as the fitting shall followed by a formula or differential equations, the finite difference method for regular grids or gradually varied functions for irregular sampling is the best way. In addition to the continuous, we can use polishing method to smooth the functions.

## 6. Discrete Smoothness, Differentiability and Lipschitz Continuity

Natural smoothness only counts the number of extreme points in general. It may not give a detailed description of discrete or digital functions. Based on the property of a polynomial, we know (n+1)-th order of derivative of a degree (n) polynomial will be zero.

A discrete function is always near a polynomial in terms of approximation as we described above. In this section, we will define a discrete smoothness to simulate the smoothness for polynomials.

Let's consider a function in the equal sampling system. Or a function $f: \{1,2,...,n\} \rightarrow R$.

$$| f(x) - f(y) | \leq Lip | x - y | \qquad (6.1)$$

Let *Lip* be the smallest constant satisfying the above equation. Let $f(x)=f^{(0)}(x)$ and $Lip=Lip^{(0)}$. So we can define the simple difference $f^{(k+1)}(x) = f^{(k)}(x+1) - f^{(k)}(x)$ and the Lipschitz constant for $f^{(k+1)}$

$$Lip^{(k)} = \max | f^{(k)}(x) - f^{(k)}(y) | / | x - y | \qquad (6.2)$$

Since if we assume that x>y,

$$| f^{(k)}(x) - f^{(k)}(y) | = | f^{(k)}(x) - f^{(k)}(x-1) + f^{(k)}(x-1) - f^{(k)}(x-2) + ... + f^{(k)}(y+1) - f^{(k)}(y) |$$
$$\leq | x - y | \max | f^{(k)}(u+1) - f^{(k)}(u) |$$

Where *(y+1)<u<x*. So, we can prove that

$$Lip^{(k)} = \max | f^{(k+1)}(x) | \qquad (6.3)$$

What we expect is that there exist *k0* such that when *k>k0*, we have

$$Lip^{(k)} > Lip^{(k+1)} \qquad (6.4)$$

**Definition 6.1** (absolute) A discrete function is called a discrete smooth function if there is a *K* such that $Lip^{(k)}$ defined in (6.1) is zero when *k>K*.

**Definition 6.2** (almost) A discrete function is called a discrete smooth function if there is a $K$ such that $Lip^{(k)}$ defined in (6.1) is smaller than $c2/(2^{(k-c1)})$ when $k>K$, $c1$, $c2$ are constants.

**Definition 6.3** (k-discrete smooth function) A discrete function is called a K-order discrete smooth function if there is a $K$ such that $Lip^{(k)}$ defined in (6.1) is smaller than $c2/(2^{(k-c1)})$ when $k<=K$, $c1$, $c2$ are constants. $c2$ is linear to $Lip^{(0)}$.

The analysis based on the above definitions will be posted next.

## References

1. L. Chen, Digital-discrete method for smooth-continuous data reconstruction, Capital Science 2010 of The Washington Academy of Sciences and its Affiliates, March 27 – 28, 2010.
2. L. Chen, Digital-Discrete Surface Reconstruction: A true universal and nonlinear method, http://arxiv.org/ftp/arxiv/papers/1003/1003.2242.pdf.
3. L. Chen, Gradual variation analysis for groundwater flow of DC (revised), DC Water Resources Research Institute Final Report, Dec 2009. http://arxiv.org/ftp/arxiv/papers/1001/1001.3190.pdf
4. L. Chen, Discrete surfaces and manifolds, Scientific and Practical Computing, Rockville, Maryland, 2004.
5. L. Chen, Gradually varied surfaces and gradually varied functions, in Chinese, 1990; in English 2005 CITR-TR 156, U of Auckland.
6. L. Chen, The necessary and sufficient condition and the efficient algorithms for gradually varied fill, Chinese Sci. Bull. 35 (10) (1990) 870-873.
7. L. Chen, Random gradually varied surface fitting, Chinese Sci. Bull. 37 (16) (1992) 1325-1329.
8. L. Chen and O. Adjei, lambda-connected segmentation and fitting, Proceedings of IEEE international conference on systems man and cybernetics, Vol 4, pp 3500-3506, 2004.
9. L. Chen, Y. Liu and F. Luo, A note on gradually varied functions and harmonic functions, 2009, http://arxiv.org/PS_cache/arxiv/pdf/0910/0910.5040v1.pdf
10. A.C. DimmiCk, Effects of sheet moisture and calender pressure on PCC and GCC coated papers, Minerals Technologies Inc. 2007.
11. H. Federer. *Geometric Measure Theory*. Springer, Berlin 1969. Page 202.
12. C. Fefferman, Whitney's extension problems and interpolation of data, Bull. Amer. Math. Soc. 46 (2009), 207-220.
13. M. D. Kirszbraun. *Über die zusammenziehende und Lipschitzsche Transformationen*. Fund. Math., (22):77–108, 1934.
14. S. G. Krantz, Lipschitz spaces, smoothness of functions, and approximation theory, Expositions Mathematicae, Vol 3, 1983, pp 193-260.
15. E. J. McShane, Extension of range of functions, Bull. Amer. Math. Soc., 40:837-842, 1934.
16. A. Rosenfeld, 'Continuous' functions on digital pictures, Pattern Recognition Letters, Vol. 4, July 1986, pp 177-184
17. J.M. Steele, Certifying smoothness of discrete functions and measuring legitimacy of images, J. of Complexity 5, 261-270 (1989).
18. J.T. Schwartz. *Nonlinear functional analysis*. Gordon and Breach Science Publishers, New York, 1969.

19. C. J. Taylor, Advanced Machine Clothing to Optimize Board Smoothness and Machine Efficiencies. In 57th Appita Annual Conference and Exhibition, Melbourne, Australia 5-7 May 2003 Proceedings. Carlton, Vic.: Appita Inc., 2003: 17-23.
20. W. Thurston, Three-dimensional geometry and topology, Princeton University press, 1997.
21. F. A. Valentine, "On the extension of a vector function so as to preserve a Lipschitz condition," Bulletin of the American Mathematical Society, vol. 49, pp. 100–108, 1943.
22. F. A. Valentine, "A Lipschitz Condition Preserving Extension for a Vector Function," American Journal of Mathematics, Vol. 67, No. 1 (Jan., 1945), pp. 83-93.
23. H. Whitney, Analytic extensions of functions defined in closed sets, *Transactions of the American Mathematical Society* **36**: 63–89, 1934.